# INFINITE-SERVER M|G|∞ QUEUEING MODELS WITH CATASTROPHES

K. Kerobyan

The infinite-server queueing models M|G|∞, BM|G|∞, $BM_k|G_k|\infty$ with homogeneous and non-homogeneous arrivals of customers and catastrophes are considered. The probability generating functions (PGF) of joint distributions of numbers of busy servers and different type served customers, as well as the Laplace-Stieltjes Transforms (LST) of distribution of busy period and distribution of busy cycle for the models are found.

**Introduction.** The infinite-server models M|G|∞ with Poisson arrival of customers and general service time distribution [1- 3] has been widely used in different fields, such as mathematical ecology, biology, finances, physics, reliability theory, industrial engineering, traffic engineering, transportation, etc, due to their following valuable properties: the stationary distribution of the model is insensitive to the form of service time distribution; the departure process of served customers follow Poisson distribution [4].

The infinite-server model M|G|∞ is one of the popular queuing models. There are numbers of generalizations of this model for different arrival processes: Batch Poisson Processes, Markov Arrival Processes (MAP), Phase Type Processes (PH), General (G) or Renewal Processes (RN) [4-17]. Transient distribution of the number of busy servers in M|M|∞ model and transient probabilities of corresponding Markov process were derived by Riordan [5]. In [6], Sevastyanov presented his famous theorem of insensitivity of stationary distribution of M|G|∞ model. He proved that stationary distribution of this model is insensitive (invariance) to the shape of service-time distribution. In [7], Benes carried out the transient probabilities of M|G|∞ model. By using probabilistic method and embedded Markov chain method, Takacs [8] derived the probability generating functions (PGF) of the number of busy servers for M|G|∞ and G|M|∞ models. Using differential equations and Laplace transformations (LT) in [9], Shanbhag derived the joint distribution of the number of busy servers, the number of served customers and total occupation time of servers at moment $t$ for the model BM|G|∞ with batch arrival of customers. In [10], Brown and Ross generalized these results by considering BM|G|∞ model with non-homogeneous Poisson arrival of customers and assumed that both the batch size and service time distributions might depend on the arrival time. By using probabilistic methods, they derived PGF of the number of busy servers and the number of served customers. In [11], Klimov considered many infinite-server models and found their main performance indices, e.g. transient probabilities of Markov process for the M|G|∞ model and the integral equations for PGF of queue of the model G|G|∞. Queue size distribution of the models G|G|∞ and BG|G|∞ with general distributions of arrival and service processes have been considered by Liu in [12]. Busy period distributions for M|G|∞ and G|G|∞ models have been derived in [13]. Model of PH|G|∞ queue with Phase Type arrival have been considered by Neuts and Ramaswami in [14]. The Basic System of Differential Equations for queue size and its matrix exponential solution are found. The differential equation for the PGF of queue size and number of served customers of the model BMAP|G|∞ with homogeneous and non-homogeneous Batch Markov Arrival of customers have been derived by Breuer [15], Kerobyan [16]. In [17] Tong found the PGF of number of busy servers for the model $M_n|G_n|\infty$ with n types of Poisson arrival of customers. This result was generalized by Masuyama and Takine in [18] for $MAP_n|G_n|\infty$ model when n-types of customers arrive according to Markov Arrival Process. By Nazarov, Moiseev and all [19-21] have been considered infinite-server queueing models with general service time distribution and different arrival processes HM|G|∞, MMP|M|∞, MAP|G|∞, and HIGI|G|∞ under condition of high intensity traffic. By using arrival process special thinning procedures, Kolmogorov differential equations, and asymptotic methods they defined first and second order approximations for PGF of queue size. However, this method cannot be applied to the models with batch arrival processes.

For many applications it is important to evaluate the inference of environmental parameters on their performance measures. For example, QoS of mobile networks has high correlation with external and internal environmental parameters. Moreover, the reliability of cables has high inference on throughput of channels and loss of requests in wired networks [22].

To evaluate the QoS of systems when the inference of environment has "disastrous" or "catastrophe" character, for example loss of communication, failure of all channels or servers of the system, instantaneous death of all species of population, collapse of financial organization or insurance company, virus attack of computer server etc., the queueing models with "catastrophes" are used. Comprehensive review on queueing models with catastrophes and negative signals can be found in Artalejo [22], Do [23], Bocharov [24]. Most of the articles are dedicated to M|M|N, G|M|1 and M|G|1 queueing models with catastrophes. Infinite-server queuing models M|M|∞ with catastrophes have been considered by Chao [25], Bohn [26], Economou and Fakinos [27], Di Crescenzo and all [28]. The stationary and transient distributions of queue size were obtained by using Markov Processes, Renewal Processes, dual processes, and embedded processes. The models of M|G|∞ and BM|G|∞ with homogeneous and non-homogeneous arrival of catastrophes are carried out in Kerobyan [29]. The basic differential equations, their solutions and LT of busy period and busy cycle of models are found. The queueing model $MMAP_n|G_n|\infty$ with homogeneous and non-homogeneous Marked MAP arrival of different types of customers and catastrophes is considered by Kerobyan [30]. The basic system of differential equations for PGF of the model, their solutions, the PGF of joint distribution of queue size and the vector of number of already solved customers are found. The queueing model $MMAP_n|G_n|\infty$ with homogeneous Marked MAP arrival of different types of customers in semi-Markov environment, with resource vector of customers and catastrophes is considered by Kerobyan [31]. The basic system of differential equations for PGF of the model, their solutions, the PGF of joint distribution of accumulated resources, queue size, vector of total served resources are found. The general method of modeling of queuenig models in random environment and catastrophes is considered.

In the present paper, we consider the models M|G|∞, BM|G|∞ and $BM_k|G_k|\infty$ with stationary and transient arrivals of customers and catastrophes. The basic differential equations of the models and their stationary and transient solutions are found. For these models the joint distribution of queue size and number of served customers and their moments are found. The distributions of busy period and busy cycles of models with stationary arrivals of customers and catastrophes have been investigated.

**2º . Model description.** Let consider an infinite-server model M|G|∞ with catastrophes where the customers and catastrophes arrive to the model according to a Poisson distribution with parameters $\lambda$ and $\nu$, respectively. The service times of the customers $\beta$ are independent identically distributed (i.i.d.) random variables (r.v.) which have a general distribution function (DF) $B(t) = P(\beta \leq t)$ with finite mean $b_1$ value. The model has infinite number of servers, and the service of arriving customer starts immediately. If a catastrophe occurs when the model is busy, then all customers in the model are immediately lost. If a catastrophe occurs when the model is empty, then it disappears without any consequences. Let the r.v. $N(t)$, $t \geq 0$ be the number of customers in the model at moment $t$, $N(t) = \{0,1,2,...\}$. Suppose that at starting time $t = 0$ the model is empty, $N(0) = 0$.

Let introduce probabilities $P_n(t) = P\{N(t) = n \mid N(0) = 0\}$, $n \geq 0$, where $P_n(t)$ is a conditional probability of having $n$, $n \geq 0$ busy servers in the model at moment $t$, if at starting time $t = 0$ the model is empty.

**Theorem 1**. The conditional probabilities of the model $P_n(t)$ satisfy following Kolmogorov differential equations

$$\frac{d}{dt}P_0(t) = -[\lambda(1-B(T-t))+v]P_0(t)+v$$

$$\frac{d}{dt}P_n(t) = -[\lambda(1-B(T-t))+v]P_n(t)+\lambda(1-B(T-t))P_{n-1}(t), \quad n \geq 1.$$

with initial conditions $P_0(0)=1$, $P_n(0)=0$, $n \geq 1$.

***Proof.*** Let consider possible changes of stochastic process $N(t)$ during the time interval $(t,t+\Delta)$. Arriving at moment $t$ customer still be in the system at moment $T$ with probability $\bar{S}(t)=1-B(T-t)$ and completed its service and left the model before moment $T$ with probability $S(t)=B(T-t)$. By the standard infinitesimal method Tijms [3] for probabilities $P_n(t)$ we derive

$$\frac{d}{dt}P_0(t) = -\lambda\bar{S}(t)P_0(t)+v\sum_{n=1}^{\infty}P_n(t)$$
$$\frac{d}{dt}P_n(t) = -[\lambda\bar{S}(t)+v]P_n(t)+\lambda\bar{S}(t)P_{n-1}(t), \quad n \geq 1 \tag{1}$$

with initial conditions $P_0(0)=1$, $P_n(0)=0$, $n \geq 1$.

The probabilities $P_n(t)$ satisfy following conditions: $0 \leq P_n(t) \leq 1$, $n \geq 0$, $t \geq 0$,

$$\sum_{n=0}^{\infty}P_n(t)=1 \tag{2}$$

From condition (2) we derive

$$1-P_0(t) = \sum_{n=1}^{\infty}P_n(t) \tag{3}$$

After substitution (3) into (1) for probabilities we get

$$\frac{d}{dt}P_0(t) = -[\lambda\bar{S}(t)+v]P_0(t)+v$$
$$\frac{d}{dt}P_n(t) = -[\lambda\bar{S}(t)+v]P_n(t)+\lambda\bar{S}(t)P_{n-1}(t), \quad n \geq 1. \tag{4}$$

with initial conditions $P_0(0)=1$, $P_n(0)=0$, $n \geq 1$.

**3º. Model Analysis**. **Transient Probabilities of the States $P_n(t)$.** To solve the system of differential equations (4), we use probability generating functions (PGF) [3,10, 35]. Let $P(z,t)$ be a PGF for number of busy servers in the model at moment $t$,

$$P(z,t) = \sum_{n=0}^{\infty}z^n P_n(t), \quad |z| \leq 1.$$

***Theorem 2***. The generating function $P(z,t)$ satisfies the following basic differential equation

$$\frac{d}{dt}P(z,t) = -[\lambda\bar{S}(t)(1-z)+v]P(z,t)+v, \quad |z| \leq 1 \tag{5}$$

with initial condition $P(z,0) = 1$.

The solution for $P(z,t)$ is given by

$$P(z,t) = e^{-\int_0^t [\lambda(1-B(x))(1-z)+v]dx} \left(1 + v\int_0^t e^{\int_0^u [\lambda(1-B(x))(1-z)+v]dx} du\right), \quad |z| \leq 1 \qquad (6)$$

or

$$P(z,t) = e^{-\int_0^t [\lambda(1-B(x))(1-z)+v]dx} + v\int_0^t e^{-\int_u^t [\lambda(1-B(x))(1-z)+v]dx} du.$$

The first $n$-th order moments $m_n(t)$ of number of busy servers can be found from (5). It is well known [3] that $m_n(t)$ satisfies

$$m_n(t) = \frac{\partial^n}{\partial z^n} P(z,t)\bigg|_{z=1}.$$

For $m_n(t)$ from (5) we get the following differential equations

$$\frac{d}{dt}m_1(t) + vm_1(t) = \lambda \overline{S}(t),$$

$$\frac{d}{dt}m_2(t) + vm_2(t) = 2\lambda \overline{S}(t)m_1(t), \qquad (7)$$

$$\frac{d}{dt}m_n(t) + vm_n(t) = n\lambda \overline{S}(t)m_{n-1}(t).$$

With initial conditions: $m_1(0) = 0$, $m_2(0) = 0$, ..., $m_n(0) = 0$.

The solutions of (7) have form

$$m_1(t) = \lambda \int_0^t \overline{S}(t) e^{-v(t-x)} dx, \qquad (8)$$

$$m_n(t) = n\lambda \int_0^t \overline{S}(t) m_{n-1}(x) e^{-v(t-x)} dx.$$

$$m_n(t) = n! \lambda^n e^{-vt} \int_0^t \overline{S}(x_n) \int_0^{x_n} \overline{S}(x_{n-1}) ... \int_0^{x_2} \overline{S}(x_1) e^{vx_1} dx_n dx_{n-1} ... dx_1.$$

Particularly, for first two moments of number of busy servers from (8) we find

$$m_1(t) = \lambda \int_0^t S(x) e^{-v(t-x)} dx, \qquad (9)$$

$$m_2(t) = 2\lambda^2 e^{-vt} \int_0^t \overline{S}(x_2) \int_0^{x_2} \overline{S}(x_1) .. e^{vx_1} dx_1 dx_2. \qquad (10)$$

The probabilities $P_n(t)$ of the model can be found by

$$P_n(t) = \frac{d^n}{dz^n} P(z,t)\Big|_{z=1}.$$

**Theorem 3.** The transient probabilities of the model $P_n(t)$ are given by

$$P_n(t) = e^{-\rho\delta(t)-vt}\frac{\rho^n\delta^n(t)}{n!} + v\frac{\rho^n}{n!}\int_0^t e^{-[\rho(\delta(t)-\delta(u))+v(t-u)]}[\delta(t)-\delta(u)]^n du. \qquad (11)$$

where $\delta(t) = \frac{1}{b_1}\int_0^t (1-B(t-x))dx = \frac{1}{b_1}\int_0^t (1-B(x))dx$, $\rho = \lambda b_1$.

Let consider some performance measures for the model M|D|∞ with catastrophes.

$$1 - B(x) = \begin{cases} 1 & \text{if } x \leq b, \\ 0 & \text{if } x > b. \end{cases}$$

$$P_n(t) = \begin{cases} \dfrac{(\lambda t)^n}{n!}e^{-(\lambda+v)t} + v\displaystyle\int_0^t e^{-[\lambda+v](t-u)}\dfrac{[\lambda(t-u)]^n}{n!}du & \text{if } t \leq b, \\[2mm] \dfrac{(\lambda b)^n}{n!}e^{-\lambda b-vt} + v\displaystyle\int_0^b e^{-[\lambda(b-u)+v(t-u)]}\dfrac{[\lambda(b-u)]^n}{n!}du & \text{if } t > b, \end{cases}$$

$$m_1(t) = \begin{cases} \dfrac{\lambda}{v}\left[1-e^{-vt}\right], & \text{if } t \leq b, \\ 0, & \text{if } t \leq b, \end{cases}$$

$$m_2(t) = \begin{cases} 2\lambda^2\left[\dfrac{(1-e^{-vt})}{v^2} - \dfrac{t}{v}e^{-vt}\right], & \text{if } t \leq b, \\ 0, & \text{if } t \leq b. \end{cases}$$

In many applications, for example in mobile computer networks and telecommunication networks, customers and catastrophes arrivals have non-stationary nature [2, 16, 34]. This fact leads to use the non-stationary Poisson processes to model the arrival of customers and catastrophes. To define the PGF $P(z,t)$, as well as the transient probabilities $P_n(t)$, first and second moments of a number of busy servers $m_1(t)$ and $m_2(t)$ at moment $t$ in the model for non-stationary arrival of catastrophes, in (6), (9), (10) and (11) we have to change $v$ into $v(t)$.

$$P(z,t) = e^{-\int_0^t [\lambda(1-B(x))(1-z)+v(x)]dx}\left(1 + \int_0^t v(x)e^{\int_0^u [\lambda(1-B(x))(1-z)+v(x)]dx}du\right), \quad |z| \leq 1. \qquad (12)$$

$$m_1(t) = \lambda\int_0^t (1-B(x))e^{-\int_x^t v(u)du}dx \qquad (13)$$

$$m_2(t) = 2\lambda\int_0^t (1-B(x))m_1(x)e^{-\int_x^t v(u)du}dx$$

$$P_n(t) = e^{-\rho\delta(t)-\int_0^t v(x)dx} \frac{\rho^n \delta^n(t)}{n!} + \frac{\rho^n}{n!}\int_0^t v(u)e^{-[\rho(\delta(t)-\delta(u))+\int_u^t v(x)dx]}[\delta(t)-\delta(u)]^n du. \tag{14}$$

When $v=0$, we have queuing model M|G|∞ with infinite number of channels and without catastrophes. From (6), (9), and (10), we derive the results that are well known in queuing theory [3, 8]

$$P(z,t) = e^{-\int_0^t [\lambda(1-B(x))(1-z)]dx}, \quad |z|\leq 1.$$

$$P(n,t) = \frac{\left(\lambda\int_0^t (1-B(x))dx\right)^n}{n!} e^{-\lambda\int_0^t (1-B(x))dx}, \quad n\geq 0.$$

$$m_1(t) = \lambda\int_0^t (1-B(x))dx, \quad m_2(t) = \left(\lambda\int_0^t (1-B(x))dx\right)^2$$

In many applications, to define the number of channels necessary for transmission of a given amount of information the distribution of the number of served customers on time interval $[0,t)$ and its moments can be used.

Let r.v. $M(t)$ be the number of customers served during time interval $[0,t)$. Define the joint distribution of r.v. $M(t)$ and $N(t)$. Let $P_{nm}(t)$ be a probability that there are $m$ busy servers at moment $t$, and $n$ customers are served on $[0,t)$ interval: $P_{mn}(t) = P(N(t)=m, M(t)=n), \; n,m \geq 0$.

Let $P(z, y, t)$ be a joint PGF of the number of busy servers at moment $t$ and number of served customers in $[0,t)$ interval.

$$P(y,z,t) = \sum_{n=0}^{\infty}\sum_{m=0}^{\infty} z^n y^m P_{mn}(t), \quad |z|\leq 1, \; |y|\leq 1.$$

**Theorem 4**. The generating function of stochastic process $(N(t), M(t))$ satisfy the following basic differential equations

$$\frac{\partial}{\partial t}P(y,z,t) = -[\lambda+v-\lambda zB(t)-\lambda y(1-B(t))]P(y,z,t) + ve^{-\lambda\int_0^t B(x)(1-z)dx}, \quad |z|\leq 1, |y|\leq 1 \tag{15}$$

with initial conditions $P(y,z,0)=1, \; |z|\leq 1, |y|\leq 1.$

**Proof.** Let consider transitions of stochastic process $(N(t), M(t))$ during the time interval $(t, t+\Delta)$. By the standard method Tijms [3] we can write the corresponding Kolmogorov differential equations:

$$\frac{d}{dt}P_{00}(t) = -(\lambda+v)P_{00}(t) + v\sum_{n=0}^{\infty} P_{n0}(t)$$

$$\frac{d}{dt}P_{0n}(t) = -(\lambda+v)P_{0n}(t) + \lambda S(t)P_{0n-1}(t) + v\sum_{m=0}^{\infty} P_{mn}(t), \quad n\geq 1$$

$$\frac{d}{dt}P_{m0}(t) = -(\lambda+v)P_{m0}(t) + \lambda\overline{S}(t)P_{m-10}(t), \quad m \geq 1$$

$$\frac{d}{dt}P_{mn}(t) = -(\lambda+v)P_{mn}(t) + \lambda S(t)P_{mn-1}(t) + \lambda\overline{S}(t)P_{m-1n}(t), \quad n \geq 1, \ m \geq 1.$$

(16)

with initial conditions $P_{00}(0) = 1, \ P_{nm}(0) = 0, \ n \geq 0, \ m \geq 0$.

Then from (16) for PGF $P(y,z,t)$ we get the following differential equation

$$\frac{\partial}{\partial t}P(y,z,t) = -[\lambda + v - \lambda zS(t) - \lambda y\overline{S}(t)]P(y,z,t) + vP(1,z,t), \quad |z| \leq 1, |y| \leq 1 \tag{17}$$

with initial conditions $P(y,z,0) = 1, \ |z| \leq 1, |y| \leq 1$.

From the following differential equation, we find the PGF $P(1,z,t)$

$$\frac{\partial}{\partial t}P(1,z,t) = -\lambda S(t)(1-z)P(1,z,t), \quad |z| \leq 1, \tag{18}$$

with initial conditions $P(z,1,0) = 1, \ |z| \leq 1$.

Note that $P(1,z,t)$ is a PGF of the number of served customers during time interval $[0,t)$ and, as follows from (18), it does not relate to the process of catastrophes.

The solution for differential equation (18) with corresponding initial conditions has a form of

$$P(1,z,t) = e^{-\lambda \int_0^t B(x)(1-z)dx}, \quad |z| \leq 1. \tag{19}$$

**Theorem 5.** The solution $P(y,z,t)$ of basic differential equations (15) is given by

$$P(y,z,t) = e^{-\int_0^t \{\lambda(1-[y(1-B(x))+zB(x)])+v\}dx} \left(1 + \int_0^t v e^{\int_0^u [\lambda(1-B(x))(1-y)+v]dx} du\right), \quad |z| \leq 1, \ |y| \leq 1. \tag{20}$$

For the model with non-stationary arrival of catastrophes the PGF $P(y,z,t)$ is given by

$$P(z,y,t) = e^{-\int_0^t \{\lambda(1-[y(1-B(x))+zB(x)])+v(x)\}dx} \left(1 + \int_0^t v(u)e^{\int_0^u [\lambda(1-B(x))(1-y)+v(x)]dx} du\right), \quad |z| \leq 1, \ |y| \leq 1. \tag{21}$$

**Corollary.** From (15) and (19), it follows that PGF $P(y,z,t)$ can be presented in a factor form

$$P(y,z,t) = P(1,z,t)P(y,1,t).$$

Here, $P(y,1,t)$ is a PGF of the number of busy servers, and $P(1,z,t)$ is a PGF of a number of customers served during $[0,t)$ interval. It can be shown that queueing model M|G|∞ with catastrophes is the only one which allows this decomposition.

From (20) when $v = 0$, for $P(y,z,t)$ we received well known results for the model M|G|∞ without catastrophes Matveyev and Ushakov [35]:

$$P(y,z,t) = e^{-\lambda \int_0^t [1-zB(x)-y(1-B(x))]dx}, \quad |y| \leq 1, |z| \leq 1.$$

**Remark 1.** The results can be generalized for the model BM(t)|G|∞ with batch arrival of customers. Suppose the batches of customers arrive according to non-homogeneous Poisson process with parameter $\lambda(t)$ [8,10]. The number of customers in the batch is independent, and identically distributed (i.i.d.) r.v. $\eta$ with a distribution $q_r(t) = P(t, \eta = r)$, $r = 1, 2, ...$ and with PGF $Q(z,t) = \sum_{r=0}^{\infty} z^r q_r(t)$. The service times of customers are i.i.d. r.v. with general distribution function $B(t)$ and finite mean $b_1$. Let consider the model with non-stationary Poisson arrival of catastrophes with parameter $v(t)$. Let $P(y,z,t)$ be the PGF of stochastic process $(N(t), M(t))$, where $N(t)$ is a number of busy servers at moment $t$ and $M(t)$ is a number of served customers in interval $[0,t)$. By using standard infinitesimal arguments we can define the PGF of joint distribution of r.v. $N(t)$ and $M(t)$.

**Theorem 6.** The PGF $P(y,z,t)$ of the model B(t)M|G|∞ is given by

$$P(y,z,t) = e^{-\int_0^t \{\lambda(x)[1-Q(y(1-B(t-x))+zB(t-x))]+v(x)\}dx} \left(1 + \int_0^t v e^{\int_0^u \{\lambda(x)[1-Q(y(1-B(u-x))+zB(u-x))]+v(x)\}dx} du\right), \quad |z| \leq 1, |y| \leq 1. \quad (22)$$

Using (22) we can define many performance metrics of the model: the probabilities $P_{mn}(t)$ of having $m$ busy servers at moment $t$ and $n$ served customers in interval $[0,t)$, the moments of busy servers and served in interval $[0,t)$ customers, and probability of idle state of the model $P_0(t)$. In particular, for the model with $v(t) = v$, $\lambda(t) = \lambda$ for probability $P_0(t)$ we get

$$P_0(t) = e^{-\int_0^t [\lambda(1-Q(B(x)))+v]dx} \left(1 + v \int_0^t e^{\int_0^u [\lambda(1-Q(B(x)))+v]dx} du\right). \quad (23)$$

**Remark 2.** Let consider the queuing model $M_k|G_k|\infty$ with $k$ types of customers and catastrophes which generalize the Tong [17] results. The arrivals of customers and catastrophes follow the Poisson distribution with parameters $\lambda_i$, $i = 1, 2, .., k$ and $v$. Service time of $i$-th type of customers is i.i.d. r.v. with general distribution $B_i(t)$, $i = 1, 2, .., k$, and finite first moment $b_{i1}$. Suppose that at its arrival epoch, the catastrophe destroys all the customers in the model if the model is busy. Let r.v. $N_i(t)$, $t \geq 0$ is a number of $i$ type customers in the model at $t$ moment, $N_i(t) = \{0, 1, 2, ..\}$, $i = 1, 2, .., k$ and r.v. $M_i(t)$ is a number of served $i$ type customers in interval $[0,t)$. At initial time $t = 0$, the model is empty, $N_i(0) = 0$, $M_i(0) = 0$, $i = 1, 2, .., k$.

Let introduce following notations and definitions: $\mathbf{z} = (z_1, z_2, ..., z_k)$, $\mathbf{y} = (y_1, y_2, ..., y_k)$,

$\mathbf{N}(t) = (N_1(t), N_2(t), ..., N_k(t))$, $\mathbf{n} = (n_1, n_2, ..., n_k)$, $\mathbf{M}(t) = (M_1(t), M_2(t), ..., M_k(t))$, $\mathbf{m} = (m_1, m_2, ..., m_k)$,

$P_{\mathbf{nm}}(t) = P\{\mathbf{N}(t) = \mathbf{n}, \mathbf{M}(t) = \mathbf{m} | \mathbf{N}(0) = 0, \mathbf{M}(0) = 0\}$, $n_i \geq 0$, $m_i \geq 0$, $i = 1, 2, .., k$, where $P_{\mathbf{nm}}(t)$ is a conditional probability of the event: in the model there are $n_1$ customers of the first type, $n_2$ customers of the second type, …, $n_k$ customers of the $k$-th type at moment $t$, and $m_1$ customers of the first type, $m_2$ customers of the second type, …, $m_k$ customers of the $k$-th type are served in interval $[0,t)$, if at initial moment $t = 0$ the model was empty and there are not any served customers.

Denote by $P(z, y, t)$ the joint PGF of number of busy servers $N(t)$ in the model at moment $t$ and number of $M(t)$ customers served in interval $[0,t)$.

$$P(z, y, t) = \sum_{n_1=0}^{\infty} \cdots \sum_{n_k=0}^{\infty} \sum_{m_1=0}^{\infty} \cdots \sum_{m_k=0}^{\infty} z_1^{n_1} z_2^{n_2} \cdots z_k^{n_k} y_1^{m_1} y_2^{m_2} \cdots y_k^{m_k} P_{nm}(t), \ |z| \leq 1, \ |y| \leq 1..$$

**Theorem 6.** The PGF $P(z, y, t)$ of the model $M_k|G_k|\infty$ is given by

$$P(z, y, t) = e^{-\int_0^t \left[ \sum_{i=1}^k \lambda_i [(1-B_i(x))z_i + B_i(x) y_i] + v \right] dx} \left( 1 + v \int_0^t e^{\int_0^u \left[ \sum_{i=1}^k \lambda_i [(1-B_i(x))z_i + B_i(x) y_i] + v \right]} du \right). \tag{24}$$

*Proof.* By using standard infinitesimal arguments for PGF of stochastic process ($N(t), M(t)$) we can write the basic differential equations

$$\frac{\partial}{\partial t} P(y, z, t) = -[\lambda + v - \lambda z S(t) - \lambda y \bar{S}(t)] P(y, z, t) + v P(1, z, t), \ |z| \leq 1, |y| \leq 1$$

with initial conditions $P(y, z, 0) = 1, \ |z| \leq 1, |y| \leq 1$.

The solution for differential equation (18) with corresponding initial conditions has a form of

$$P(1, z, t) = e^{-\lambda \int_0^t B(x)(1-z) dx}, \ |z| \leq 1.$$

From (24) we can define the probability $P_{nm}(t)$, the moments of the numbers of each type customers in the model at moment $t$, numbers of each type customers served in interval $[0,t)$, and probability of idle state of the model $P_0(t)$. In particular, for $P_0(t)$ we find

$$P_0(t) = e^{-\int_0^t \left[ \sum_{i=1}^k \lambda_i (1-B_i(x)) + v \right] dx} \left( 1 + v \int_0^t e^{\int_0^u \sum_{i=1}^k \lambda_i (1-B_i(x)) dx} du \right). \tag{25}$$

**Remark 3.** Now let consider the queuing model $BM_k|G_k|\infty$ with batch arrival of different types of customers and catastrophes. The arrivals of catastrophes and customers follow a Poisson distribution with parameters $v$ and $\lambda_i(n), \ i = 1, 2, .., k, n = 0, 1, 2, ....$, where $n$ is the number of customers in the batch, $\lambda_i(n) = \lambda_i P(\eta_i = n)$. The number of $i$-th type customers in the batches is i.i.d. r.v. $\eta_i$ with a distribution $q_{ir} = P(\eta_i = r), \ r = 1, 2, ...$ and with PGF $Q_i(z_i) = \sum_{r=0}^{\infty} z_i^r q_{ir}$.

Service time of $i$-th type customers is i.i.d. r.v. with general distribution $B_i(t), \ i = 1, 2, .., k$, and finite first moment $b_{i1}$. Suppose that at its arrival epoch, the catastrophe destroys all the customers in the model. Let r.v. $N_i(t), \ t \geq 0$ is a number of $i$-th type customers in the model at $t$ moment, $N_i(t) = \{0, 1, 2, ..\}, \ i = 1, 2, ..., k$. At initial time $t = 0$, the model is empty, $N_i(0) = 0, \ i = 1, 2, .., k$. For this model the PGF of joint distribution of r.v. $N(t)$ and $M(t)$ $P(y, z, t) = P(y_1, y_2, ..., y_k, z_1, z_2, ..., z_k, t)$ can be found

$$P(y,z,t) = e^{-\int_0^t \{\sum_{i=1}^k \lambda_i[1-Q_i(y_i(1-B_i(x))+z_iB_i(x))]+v\}dx} \left(1+\int_0^t ve^{\int_0^u \{\sum_{i=1}^k \lambda_i[1-Q_i(y_i(1-B_i(x))+z_iB_i(x))]+v\}dx} du\right).$$

From PGF $P(y,z,t)$ can be defined the PGFs of numbers of served customers and numbers of busy servers at moment $t$

$$P(y,1,t) = e^{-\int_0^t \{\sum_{i=1}^k \lambda_i[1-Q_i(y_i(1-B_i(x))+B_i(x))]+v\}dx} \left(1+\int_0^t ve^{\int_0^u \{\sum_{i=1}^k \lambda_i[1-Q_i(y_i(1-B_i(x))+B_i(x))]+v\}dx} du\right).$$

$$P(1,z,t) = e^{-\int_0^t \{\sum_{i=1}^k \lambda_i[1-Q_i((1-B_i(x))+z_iB_i(x))]+v\}dx} \left(1+\int_0^t ve^{\int_0^u \{\sum_{i=1}^k \lambda_i[1-Q_i((1-B_i(x))+z_iB_i(x))]+v\}dx} du\right).$$

For this model, we can define also non-stationary arrivals of different types of customers and catastrophes. The corresponding PGF of joint distribution of numbers of served customers and numbers of busy servers at moment $t$ has form

$$P(y,z,t) = e^{-\int_0^t \{\sum_{i=1}^k \lambda_i(x)[1-Q_i(y_i(1-B_i(t-x))+z_iB_i(t-x))]+v(x)\}dx} \left(1+\int_0^t v(u)e^{\int_0^u \{\sum_{i=1}^k \lambda_i(x)[1-Q_i(y_i(1-B_i(u-x))+z_iB_i(u-x))]+v(x)\}dx} du\right).$$

Let now suppose that arriving batches can contain different types of customers. The models without catastrophes have been considered by number of authors, e.g. Fakinos [32], Tang [17], Choi and Park [33]. Let random variable $x_i$ be the number of $i$-th type customers in the batch. The vector $\boldsymbol{x} = (x_1, x_2, ..., x_k)$ is denoted the batch size and has a PGF

$$Q(\boldsymbol{x}) = Q(x_1, x_2, ..., x_k) = \sum_{n_1 \geq 0} \sum_{n_2 \geq 0} ... \sum_{n_k \geq 0} q(x_1 = n_1, x_2 = n_2, ..., x_k = n_k) z_1^{n_1} z_2^{n_2} ... z_k^{n_k}$$

Where $q(x_1 = n_1, x_2 = n_2, ..., x_k = n_k)$ is a probability that arriving batch contains $n_1$ customers of first type, $n_2$ customers of second type, … $n_k$ customers of $k$ type. The service times of customers are i.i.d. r.v. and for $i$ type of customers have distribution $B_i(x)$. Customer of type $i$ arriving moment $s$ still be in the system at moment $t$ with probability $1 - B_i(t-s)$ and quit the model after completing the service before moment $t$ with probability $B_i(t-s)$. At initial moment the model is empty.

**Theorem 7**. The PGF of the model is given by

$$P_1(z,t) = \sum_{n_1 \geq 0} \sum_{n_2 \geq 0} ... \sum_{n_k \geq 0} P(N_1(t) = n_1, N_2(t) = n_2, ..., N_k(t) = n_k) z_1^{n_1} z_2^{n_2} ... z_k^{n_k} = e^{\lambda \int_0^t [1-Q(w(\boldsymbol{x},t-s))]ds} \quad (26)$$

$$P_2(z,t) = \sum_{n_1 \geq 0} \sum_{n_2 \geq 0} ... \sum_{n_k \geq 0} P(N_1(t) = n_1, N_2(t) = n_2, ..., N_k(t) = n_k) z_1^{n_1} z_2^{n_2} ... z_k^{n_k} = e^{-\lambda \int_0^t Q(w(\boldsymbol{x},t-s))ds}$$

To solve this model we will use method of collective marks. Let us mark - color each $i$ type customer in arriving batch independently from other customers in the batch and in the system with "red" color by probability $z_i$ or "blue" color by probability $1-z_i$. Then $Q(x)$ can be interpreted as probability of event "arriving batch does not contain "blue" customers". $w_i(z_i, t-s) = B_i(t-s) + (1 - B_i(t-s))z_i$ - is a probability of event "type $i$ customer arriving at the moment $s$ still is in the system at moment $t$ and has "red" color, $w(x, t-s) = (w_1(z_1, t-s), w_2(z_2, t-s), ..., w_k(z_k, t-s))$. Then the probability of event " arriving at the moment $s$ batch $x = (n_1, n_2, ..., n_k)$ which customers still are in the system at moment $t$ does not contain "blue" color customers" is $\prod_{i=1}^{k} [B_i(t-s) + (1 - B_i(t-s))z_i]^{n_i} = \prod_{i=1}^{k} w_i(z_i, t-s)^{n_i}$. Probability of event "arriving at moment $s$ batch which customers still are in the system at moment $t$ does not contain "blue" color customers" is $R(x, s, t) = \sum_{n_1 \geq 0} \sum_{n_2 \geq 0} ... \sum_{n_k \geq 0} q(n_1, n_2, ..., n_k) \prod_{i=1}^{k} [B_i(t-s) + (1 - B_i(t-s))z_i]^{n_i} = Q(w(x, t-s))$. The batches arrive according to Poisson process with parameter $\lambda$. Let each arriving batch mark with probability $Q(w(x, t-s))$ and unmark with probability $1 - Q(w(x, t-s))$. Then, as shown by Reichelt [34], the splitting processes are non-homogeneous Poisson processes with parameters $\lambda \int_0^t Q(w(x, t-s))ds$ and $\lambda \int_0^t [1 - Q(w(x, t-s))]ds$ respectively. Hence, the probability of event "no blue customers in the system at moment $t$" is

$$P_1(z, t) = \sum_{n_1 \geq 0} \sum_{n_2 \geq 0} ... \sum_{n_k \geq 0} P(N_1(t) = n_1, N_2(t) = n_2, ..., N_k(t) = n_k) z_1^{n_1} z_2^{n_2} ... z_k^{n_k} = e^{\lambda \int_0^t [1 - Q(w(x, t-s))]ds},$$

the probability of event "no red customers in the system at moment $t$" is

$$P_2(z, t) = \sum_{n_1 \geq 0} \sum_{n_2 \geq 0} ... \sum_{n_k \geq 0} P(N_1(t) = n_1, N_2(t) = n_2, ..., N_k(t) = n_k) z_1^{n_1} z_2^{n_2} ... z_k^{n_k} = e^{-\lambda \int_0^t Q(w(x, t-s))ds}.$$

Where $P_2(z, t)$ is a PGF of the number of customers that have been served before time $t$, and $P_1(z, t)$ is a PGF of the number of customers that served in the system at the moment $t$.

The corresponding PGF for the model with catastrophes can be found by following arguments [31]. Let us note that the dynamic of the model with catastrophes between two consecutive catastrophes is the same as the dynamic of the model without catastrophes. After the catastrophes the model jumps to the idle state and continues its dynamic from that state as the model without catastrophes. If the catastrophes occur according to Poisson process with parameter $v$ then the PGF of the model with catastrophes $\hat{P}(z, t)$ we define by using collective marks method.

$\hat{P}(z, t)$ is the probability that no "blue" customers in the model with catastrophes at moment $t$. Indeed, this event can happened if catastrophes do not occur in [0,t) and no "blue" customers arrive (the probability of this event is $P(z, t)e^{-vt}$), or a catastrophe occurs at the moment $x$, $x \in [0, t)$, (the probability is $vdx$), the model jumps into idle state and during the time $t - x$ catastrophes do not occur and no "blue" customers arrive (the probability of this event is $\int_0^t P(z, t-x)e^{-v(t-x)}vdx$). By using the total probability rule for PGF $\hat{P}(z, t)$ we get

$$\hat{P}(z, t) = P(z, t)e^{-vt} + v\int_0^t P(z, t-x)e^{-v(t-x)}dx. \tag{27}$$

Let $\tilde{\hat{P}}(z,s)$ is the Laplace Transformation (LT) of PGF $\hat{P}(z,t)$

$$\tilde{\hat{P}}(z,s) = \int_0^\infty e^{-st}\hat{P}(z,t)dt.$$

Then from (19) we find

$$\tilde{\hat{P}}(z,s) = \tilde{P}(z,s+v) + \frac{v}{s}\tilde{P}(z,s+v) = \tilde{P}(z,s+v)(1+\frac{v}{s}). \tag{28}$$

If this result write in the form

$$s\tilde{\hat{P}}(z,s) = \tilde{P}(z,s+v)(s+v) \tag{29}$$

then according to Kerobyan [31] and Matveyev [35] it has simple probabilistic interpretation. Let suppose that independently of the model can occur event A according to Poisson process with parameter $s$. Then the left and right sides of (21) can be interpreted as probabilities of following events:

*left side*; "the first event A occurred when in the model with catastrophes are not any "blue" customers",

(with probability $s\tilde{\hat{P}}(z,s) = s\int_0^\infty e^{-st}\hat{P}(z,t)dt$)

*right side*; "the first event of total flow of catastrophes and event A occurred when in the model without catastrophes are not any "blue" customers".

(with probability $(s+v)\tilde{P}(z,s+v) = s\int_0^\infty e^{-(s+v)t}P(z,t)dt$).

This model generalizes the results of Fakinos [32], Tang [17], Choi and Park [33]. The moments of the number of different types of customers can be found by the standard way, see for example Fakinos [32].

**4º. Model Analysis: Busy Period.** To define the DF of busy period of the model with catastrophes, we shall use the method of collective marks [11, 17, 35, 36]. Let $P_0(t)$ be the probability that the model is empty at moment $t$. At initial time $t=0$ the model is empty, $P_0(0)=1$, and $\tilde{P}_0(s)$ is a Laplace Transform (LT) of a function $P_0(t)$

$$\tilde{P}_0(s) = \int_0^\infty e^{-st} P_0(t)dt,$$

Let r.v. $\pi$ be the length of busy period of the model. Busy period begins when the first customer arrives at the empty model and ends when the model is free of customers. The r.v. $\omega$ be the length of busy cycle and is defined as the sum of busy period $\pi$ and idle period $\varphi$ of the model. For infinite-server M|G|$\infty$ model $\varphi$ is a r.v. with exponential distribution and parameter $\lambda$. The Laplace-Stieltjes Transforms (LST) of busy period $\hat{\pi}(s)$ and busy cycle $\hat{\omega}(s)$ define

$$\hat{\pi}(s) = \int_0^\infty e^{-st} d\pi(t), \quad \hat{\omega}(s) = \int_0^\infty e^{-st} d\omega(t).$$

According to [35] for the class of conservative queuing models with Poisson arrival of customers, the probability that the model is empty, the distributions of busy period and the busy cycle can be defined from the following equations

$$s\tilde{P}_0(s) = \frac{s}{\lambda+s} + \frac{\lambda}{\lambda+s}\hat{\pi}(s)s\tilde{P}_0(s),$$

$$\hat{\omega}(s) = \frac{\lambda}{\lambda+s}\hat{\pi}(s). \tag{30}$$

Where

$$\tilde{P}_0(s) = \int_0^\infty P_0(t)e^{-st}dt, \quad P_0(t) = e^{-\int_0^t [\lambda(1-B(x))+v]dx}\left(1 + v\int_0^t e^{\int_0^u [\lambda(1-B(x))+v]dx}du\right). \tag{31}$$

Note, that the considering infinite-server model with catastrophes is conservative. Therefore, solving (30) for busy period and busy cycle distributions, we derive

$$\hat{\pi}(s) = 1 + \frac{s}{\lambda} - \frac{1}{\lambda\tilde{P}_0(s)} = 1 + \frac{s}{\lambda}(1 - \frac{1}{s\tilde{P}_0(s)}),$$

$$\hat{\omega}(s) = 1 - \frac{\lambda+s}{\tilde{P}_0(s)}. \tag{32}$$

If idle state probability of the model M|G|∞ without catastrophes note by $p_0(t)$

$$p_0(t) = e^{-\int_0^t \lambda(1-B(t-x))dx}$$

Then from (31) we find

$$\tilde{P}_0(s) = \int_0^\infty e^{-st} e^{-\int_0^t [\lambda(1-B(t-x))+v]dx}dt + v\int_0^\infty e^{-st}\int_0^t e^{-\int_u^t \lambda(1-B(t-x))dx + v(t-u)}du$$

$$= \int_0^\infty e^{-(v+s)t}p_0(t)dt + v\int_0^\infty e^{-st}\int_0^t p_0(t-u)e^{-v(t-u)}du \tag{33}$$

$$= \tilde{p}_0(v+s)(1 + \frac{v}{s}).$$

This result can be written in the form $s\tilde{P}_0 = (v+s)\tilde{p}_0(v+s)$ and interpreted by similar to (29) way.

In particular, when $v = 0$, from (30) we derive the well known result for classical model M|G|∞ without catastrophes Stadje [37]:

$$\hat{\pi}(s) = 1 + \lambda^{-1}[s - \frac{1}{\int_0^\infty e^{-\int_0^t [\lambda(1-B(x))+s]dx}dt}], \quad \hat{\omega}(s) = 1 - \frac{\lambda+s}{\int_0^\infty e^{-\int_0^t [\lambda(1-B(x))+s]dx}dt}.$$

Let consider the busy period $\hat{\pi}(s)$ and busy cycle $\hat{\omega}(s)$ distributions for the model M|D|∞ with catastrophes.

$$B(x) = \begin{cases} 0 & \text{if } x \leq b, \\ 1 & \text{if } x > b. \end{cases}$$

Then for busy period $\hat{\pi}(s)$ and busy cycle $\hat{\omega}(s)$ we found

$$\tilde{P}_0(s) = \frac{v + (\lambda + s)e^{(\lambda+s+v)b}}{s(\lambda + s + v)e^{(\lambda+s+v)b}}, \quad \lim_{s \to 0} s\tilde{P}_0(s) = \frac{v + \lambda e^{-(\lambda+v)b}}{v + \lambda}, \tag{34}$$

$$\hat{\pi}(s) = \frac{\lambda + s + ve^{(\lambda+s+v)b}}{\lambda + (s+v)e^{(\lambda+s+v)b}},$$

$$\hat{\omega}(s) = 1 - \frac{s(\lambda + s)(\lambda + s + v)e^{(\lambda+s+v)b}}{\lambda + (v+s)e^{(\lambda+s+v)b}}. \tag{35}$$

$$\bar{\pi}_1 = \frac{e^{(\lambda+v)b} - 1}{\lambda + ve^{(\lambda+v)b}},$$

$$\bar{\pi}_2 = \frac{2[1 - (1 + b(\lambda+v))e^{(\lambda+v)b}]}{[v + \lambda e^{(\lambda+v)b}]^2}. \tag{36}$$

If $v = 0$, from (19) we derive the well known result for classical model $M|D|\infty$ without catastrophes Nazarov [38]:

$$\hat{\pi}(s) = \frac{\lambda + s}{\lambda + se^{(\lambda+s)b}},$$

$$\hat{\omega}(s) = 1 - \frac{(\lambda + s)^2 s}{s + \lambda e^{-(\lambda+s)b}}.$$

The busy period and busy cycle of the models $BM|G|\infty$ and $M_k|G_k|\infty$ can be found by using the idle state probabilities of corresponding models (23) and (25). For the LST of busy period of the models $BM|G|\infty$ and $M_k|G_k|\infty$ we get

$$\hat{\pi}_{BM|G|\infty}(s) = 1 + \lambda^{-1}\left[s - \frac{1}{\int_0^\infty e^{-\int_0^t [\lambda(1-Q(B(x)))+v+s]dx}\left(1 + v\int_0^t e^{\int_0^u [\lambda(1-Q(B(x)))+v]dx} du\right)dt}\right], \tag{37}$$

$$\hat{\pi}_{Mk|Gk|\infty}(s) = 1 + \lambda^{-1}\left[s - \frac{1}{\int_0^\infty e^{-\int_0^t \left[\sum_{i=1}^k \lambda_i(1-B_i(x))+v+s\right]dx}\left(1 + v\int_0^t e^{\int_0^u \sum_{i=1}^k \lambda_i(1-B_i(x))dx} du\right)dt}\right]. \tag{38}$$

The results of this paper can be generalized in following directions by considering models with: Batch Markov Arrival Process, Marked Markov Arrival Process or more general arrival processes of customers and catastrophes; stationary and transient arrival processes of customers and catastrophes; service time distribution of customers depend on their type, time moment of their arrival and an additional parameters vector; the models can be considered in some random environment.